\documentclass[a4paper,reqno]{amsart}
\usepackage{amssymb, amsmath, amscd}

\def\gronko{\vphantom{\vrule height 12pt}}
\def\II((#1)){{\mathcal{I}}\{(#1)\}}\def\trace{\operatorname{Trace}}
\newtheorem{theorem}{Theorem}[section]
\newtheorem{definition}[theorem]{Definition}
\newtheorem{lemma}[theorem]{Lemma}

\newtheorem{remark}[theorem]{Remark}

\def\id{\operatorname{id}}

\makeatletter
 \@addtoreset{equation}{section}
\makeatother
\begin{document}
\title[Almost pseudo-Hermitian and almost para-Hermitian manifolds]{
Geometric realizability of covariant derivative K\"ahler tensors for almost pseudo-Hermitian and almost para-Hermitian
manifolds}
\author[M. Brozos-V\'{a}zquez et al.]
{M. Brozos-V\'{a}zquez \and E. Garc\'{\i}a-R\'{\i}o \and P. Gilkey \and L. Hervella}
\address{MB-V: Department of Mathematics, University of A Coru\~na, Spain}
\email{miguel.brozos.vazquez@udc.es}
\address{EG-R: Faculty of Mathematics, University of Santiago de Compostela, Spain}
\email{eduardo.garcia.rio@usc.es}
\address{PG: Mathematics Department, University of Oregon,
  Eugene OR 97403, USA}
  \email{gilkey@uoregon.edu}
\address{LH: Faculty of Mathematics, University of Santiago de Compostela, Spain}
\email{luismaria.hervella@usc.es}
\begin{abstract}{The covariant derivative of the K\"ahler form of an almost pseudo-Hermitian
or of an almost para-Hermitian manifold satisfies certain algebraic relations. We show, conversely,
that any $3$-tensor which satisfies these algebraic relations can be realized geometrically.
}\end{abstract}
\maketitle

\section{\bf Introduction}\label{sect-1}

The paper of Gray and Hervella \cite{GH80} puts into a unified framework 16 classes of almost Hermitian manifolds and was the work  which inspired other classification results like those in
\cite{Na83,TV81,TV83}.
It is important in the mathematical setting and is used in obvious settings when
some class of K\"ahler  or Hermitian  manifolds is the central focus of investigation.
The Gray-Hervella decomposition plays a role in the discussion of nearly K\"ahler and almost K\"ahler geometry as well as in the study of conformal equivalences among almost Hermitian structures (see for example \cite{EPS09,MO08}, \cite{Arm02},  and \cite{BLS07,CI08},  respectively).
It is related to the Tricerri-Vanhecke \cite{TV81} decomposition of
the curvature tensor in \cite{FFS94} and it has a prominent role in understanding the influence of the curvature on the underlying structure of the manifold \cite{HL00}.
The Gray-Hervella classification is
related to the 64 classes of almost quaternion-Hermitian structures in \cite{MS04}, showing some interactions amongst them.
The different classes have been considered for flag manifolds -- they essentially reduce to four classes \cite{MN03}, and the 6-dimensional case has been considered in detail in \cite{AFS05}. The different classes of almost Hermitian structures also enter into the discussion of some harmonicity
problems~\cite{BLS07}.

Although most of this work has been in the positive definite setting, the indefinite case
also plays a role (see for example  \cite{CS07,GM00,GK05,MHL07,SY08}). In addition to the pseudo-Hermitian setting, the almost para-Hermitian geometry is of interest both from the mathematical and the physical point of view \cite{AGK11,AMT09,CMMS04,CMMS05,GPSR07,P08}.
Related work of Gadea and Masque \cite{GM91} classified almost para-Hermitian structures into 32 different classes by considering separately the two natural distributions associated to the almost para-Hermitian structure.

In this paper we put both the almost para-Hermitian and the almost pseudo-Hermitian structures in an unified context by extending the Gray-Hervella decomposition to the pseudo-Riemannian setting. This is done by analyzing the covariant derivative of the corresponding K\"ahler form and the decomposition of the space of such tensors under the action of a suitable structure group (see Theorem
\ref{thm-1.4} for details). Moreover we consider the geometric realizability of all the different classes by perturbing the given structures.
In Theorem
\ref{thm-1.1}, we show that {\it any} algebraic covariant derivative K\"ahler tensor can be geometrically realized by perturbing the underlying structure on a given almost para/pseudo-Hermitian background manifold; Theorem \ref{thm-1.2} provides a similar result in the integrable setting. In Theorem \ref{thm-1.6}, we restrict to the complex setting and extend results of \cite{GH80} from the positive definite context to the indefinite context showing any of the 16 classes has at least one
geometrical representative.

\medskip

We establish notation as follows. Let $(M,g)$ be a pseudo-Riemannian manifold of dimension $m=2\bar m$. Let $J_\pm$ be
endomorphisms of the tangent bundle $TM$.
We say that $(M,g,J_+)$ is an {\it almost para-Hermitian
manifold} if $J_+^2=\id$ and if $J_+^*g=-g$. Similarly, if $J_-^2=-\id$ and if
$J_-^*g=g$, then we say that
$(M,g,J_-)$ is an {\it almost pseudo-Hermitian manifold}. The existence of such
structures is related to the signature
$(p,q)$ of $g$. If
$(M,g)$ admits an almost para-Hermitian structure $J_+$, then
$p=q$. Similarly if $(M,g)$ admits an almost pseudo-Hermitian structure $J_-$, then both $p$ and
$q$ are even. Thus usually we are not dealing with both $J_-$ and $J_+$ at the same time on $(M,g)$,
but we adopt a common notation to keep the exposition in parallel as much as possible.

Let $\nabla$ be the Levi-Civita connection of $g$. The {\it
associated K\"ahler form} and the covariant derivative are defined, respectively, by:
\begin{eqnarray*}
&&\Omega_\pm(x,y):=g(x,J_\pm y),\\
&&\nabla\Omega_\pm(x,y;z)=zg(x,J_\pm y)-g(\nabla_zx,J_\pm y)-g(x,J_\pm\nabla_zy)\,.
\end{eqnarray*}
We subscript $J$ and $\Omega$ to keep track of the signs involved. For example, as we
shall see presently in Lemma
\ref{lem-3.1}, we have:
\begin{equation}\label{eqn-1.a}
\begin{array}{l}
\nabla\Omega_\pm(x,y;z)=-\nabla\Omega_\pm(y,x;z),\\
\nabla\Omega_\pm(x,y;z)=\pm\nabla\Omega_\pm(J_\pm x,J_\pm y;z)\,.\gronko
\end{array}\end{equation}

It is convenient to work in an algebraic context as well. Let $(V,\langle\cdot,\cdot\rangle)$ be an inner product space and
let $J_\pm^0$ be linear maps of $V$. We say that
$(V,\langle\cdot,\cdot\rangle,J_+^0)$ is a {\it para-Hermitian
vector space} if
$(J_+^0)^*\langle\cdot,\cdot\rangle=-\langle\cdot,\cdot\rangle$ and if
$(J_+^0)^2=\id$. Similarly, $(V,\langle\cdot,\cdot\rangle,J_-^0)$
is said to be a {\it pseudo-Hermitian vector space} if
$(J_-^0)^*\langle\cdot,\cdot\rangle=\langle\cdot,\cdot\rangle$ and if
$(J_-^0)^2=-\id$. Again, the existence of such structures imposes restrictions on the signature. Motivated by
Equation~(\ref{eqn-1.a}), we define:
$$
\begin{array}{l}
\mathfrak{H}_\pm:=
 \{H_\pm\in\otimes^3V^*:H_\pm(x,y;z)=-H_\pm(y,x;z)\quad\text{and}\\
\gronko\qquad\qquad\quad H_\pm(J_\pm^0x,J_\pm^0y;z)=\pm
H_\pm(x,y;z)\ \forall\ x,y,z\}\,.
\end{array}$$
Let $H_\pm\in\mathfrak{H}_\pm$. We have
\begin{equation}\label{eqn-1.b}
H_\pm(x,J_\pm^0y;z)=\pm H_\pm(J_\pm^0x,J_\pm^0J_\pm^0y;z)=H_\pm(J_\pm^0x,y;z)\,.
\end{equation}

The following result shows that Equation~(\ref{eqn-1.a}) generates the universal
symmetries satisfied by
$\nabla\Omega_\pm$ and provides a rich family of examples. It is striking that we can fix the metric and only vary
the almost (para)-complex structure; in particular, we could take the background structure to be flat.

\begin{theorem}\label{thm-1.1}
Let $(M,g,J_\pm)$ be a background almost para/pseudo-Hermitian ma\-nifold and let $P\in M$. Suppose given
$H_\pm$ in $\mathfrak{H}_\pm(T_PM,g_P,J_{\pm,P})$. Then there exists a new almost para/pseudo-Hermitian structure
$\tilde J_\pm$ on $M$ which agrees with $J_\pm$ at $P$ so that
$\nabla\Omega_\pm(M,g,\tilde J_{\pm})(P)=H_\pm$.
\end{theorem}

We consider the following subspace:
$$U_{3,\pm}:=\{H_\pm\in\mathfrak{H}_\pm:H_\pm(x,y;z)=\mp
H_\pm(x,J_\pm^0y;J_\pm^0z)\ \forall\ x,y,z\}\,.$$
If $(M,g,J_\pm)$ is a para/pseudo-Hermitian manifold (i.e. $J_\pm$ is integrable), then $\nabla\Omega_\pm\in U_{3,\pm}$ as we
shall see presently in Lemma \ref{lem-3.2}. Conversely:
\begin{theorem}\label{thm-1.2}
Let $(M,g,J_\pm)$ be a background para/pseudo-Hermitian manifold and let $P\in M$. Suppose given $H_\pm$ in
$U_{3,\pm}(T_PM,g_P,J_{\pm,P})$. Then there exists a new para/pseudo-Hermitian metric $\tilde g$ on
$M$ which agrees with $g$ at $P$ so that
$\nabla\Omega_\pm(M,\tilde g,J_\pm)(P)=H_\pm$.
\end{theorem}

Theorems \ref{thm-1.1} and \ref{thm-1.2} are global results; it is necessary to have a starting background structure as not
every manifold admits a para/pseudo-Hermitian structure of a given signature; in general, there are topological restrictions on
$M$ for the existence of a (para)-complex structure or for the existence of a metric of signature $(p,q)$. These Theorems give
results in the category of compact manifolds. However it is a direct consequence of the Theorems that one can also restrict
attention to an open coordinate chart to get purely local results.

These results are based on a decomposition of $\mathfrak{H}_\pm$ which extends the decomposition given in \cite{GH80} in the
positive definite context. Adopt the {\it Einstein convention} and sum over repeated indices.
\begin{definition}
\rm Let $(V,\langle\cdot,\cdot\rangle,J_\pm^0)$ be a para/pseudo-Hermitian vector space.
Let $\varepsilon_{ij}:=\langle e_i,e_j\rangle$ for some basis $\{e_i\}$ for $V$. Let $\phi\in V^*$.  Let $H\in\otimes^3V^*$.
Let $\operatorname{GL}$ be the general linear group. Set:
\begin{enumerate}
\item $(\tau_1H)(x):=\varepsilon^{ij}H(x,e_i;e_j)$.
\smallbreak\item
$\sigma_\pm(\phi)(x,y;z):=\phi({J_\pm^0} x)\langle y,z\rangle-\phi({J_\pm^0} y)\langle x,z\rangle
+\phi(x)\langle{J_\pm^0} y,z\rangle-\phi(y)\langle{J_\pm^0} x,z\rangle$.
\smallbreak\item $W_{1,\pm}:=\{H_\pm\in\mathfrak{H}_\pm:H_\pm(x,y;z)+H_\pm(x,z;y)=0$ $\forall$ $x,y,z\}$.
\smallbreak\item $W_{2,\pm}:=\{H_\pm\in\mathfrak{H}_\pm:
   H_\pm(x,y;z)+H_\pm(y,z;x)+H_\pm(z,x;y)=0$ $\forall$ $x,y,z\}$.
\smallbreak\item $W_{3,\pm}:=U_{3,\pm}\cap\ker(\tau_1)$.
\smallbreak\item $W_{4,\pm}:=\operatorname{Range}(\sigma_\pm)$.
\smallbreak\item $\mathcal{O}:=\{T\in\operatorname{GL}:T^*\langle\cdot,\cdot\rangle=\langle\cdot,\cdot\rangle\}$.
\smallbreak\item $\mathcal{U}_\pm:=\{T\in\mathcal{O}:TJ_\pm^0=J_\pm^0T\}$.
\smallbreak\item $\mathcal{U}_\pm^\star:=\{T\in\mathcal{O}:TJ_\pm^0=TJ_\pm^0\text{ or }TJ_\pm^0=-J_\pm^0T\}$.
\smallbreak\item $\operatorname{GL}_\pm:=\{T\in\operatorname{GL}:TJ_\pm^0=J_\pm^0T\}$.
\smallbreak\item $\chi(T):=+1$ if $T\in\mathcal{U}_\pm$ and $\chi(T):=-1$ if $T\in\mathcal{U}_\pm^\star-\mathcal{U}_\pm$.
\end{enumerate}\end{definition}

\begin{theorem}\label{thm-1.4}
Let $m\ge6$. We have a direct sum orthogonal decomposition of $\mathfrak{H}_\pm$ and of $U_{3,\pm}$ into irreducible
inequivalent
$\mathcal{U}_\pm^\star$ modules
 in the form:
$$
\mathfrak{H}_\pm=W_{1,\pm}\oplus W_{2,\pm}\oplus W_{3,\pm}\oplus W_{4,\pm}\quad\text{and}\quad
U_{3,\pm}=W_{3,\pm}\oplus W_{4,\pm}\,.
$$
\end{theorem}

One obtains the corresponding decompositions if $m=4$ by setting $W_{1,\pm}=0$ and $W_{3,\pm}=0$.
The modules $W_{i,-}$ are also irreducible $\mathcal{U}_-$ modules so the decomposition of \cite{GH80} of $\mathfrak{H}_-$ as
a $\mathcal{U}_-$ module extends without change from the positive definite to the indefinite setting; we omit the additional
analysis this requires in the interests of brevity. The modules
$W_{i,+}$ are not, however, irreducible
$\mathcal{U}_+$ modules and thus the classification of
\cite{GM91} is a more refined one than we consider here as there are 8 factors in the decomposition rather than 4. By using
the structure group
$\mathcal{U}_+^\star$ instead of
$\mathcal{U}_+$, we shall bypass some of the technical difficulties encountered in
\cite{GM91} and this structure group is sufficient for our purposes.

The focus of Theorem \ref{thm-1.1} and of Theorem \ref{thm-1.2} is to show that {\it every} element of $\mathfrak{H}_\pm$ and
of $U_{3,\pm}$ is geometrically realizable in an appropriate context. One can, however, focus instead on the precise nature of
the classes involved. We now restrict to the complex setting. Let
$\xi$ be a $\mathcal{U}_-^\star$ submodule of
$\mathfrak{H}_-$. We say that $(M,g,J_-)$ is a {\it $\xi$-manifold} if $\nabla\Omega_-$ belongs to $\xi$ for
every point of the manifold and if $\xi$ is minimal with this property. This gives rise to the celebrated 16 classes of almost
Hermitian manifolds (in the positive definite setting) \cite{GH80}:
\begin{theorem}\label{thm-1.5}
Let $\xi$ be a submodule of
$\mathfrak{H}_-$. Then there exists an almost Hermitian $\xi$-manifold.
\end{theorem}

 We can generalize this to the indefinite setting; we shall suppose $m\ge10$
to simplify the discussion:

\begin{theorem}\label{thm-1.6}
Suppose given $(2\bar p,2\bar q)$ with $2\bar p+2\bar q\ge10$. Let $\xi$ be a submodule of
$\mathfrak{H}_-$. Then there exists a $\xi$-manifold of signature $(2\bar p,2\bar q)$.
\end{theorem}

Many of these classes have geometrical meanings which have been extensively investigated. For example:
\begin{enumerate}
\item $\xi=\{0\}$ defines the class of K\"ahler manifolds.
\item $\xi=W_{1,-}$ defines the class of nearly K\"ahler manifolds.
\item $\xi=W_{2,-}$ defines the class of almost K\"ahler manifolds.
\item $\xi=W_{3,-}$ defines the class of Hermitian semi-K\"ahler manifolds.
\item $\xi=W_{1,-}\oplus W_{2,-}$ defines the class of quasi-K\"ahler manifolds.
\item $\xi=W_{3,-}\oplus W_{4,-}=U_{3,-}$ defines the class of pseudo-Hermitian manifolds.
\item $\xi=W_{1,-}\oplus W_{2,-}\oplus W_{3,-}$ defines the class of semi-K\"ahler manifolds.
\item $\xi=\mathfrak{H}_-$ defines the class of almost pseudo-Hermitian manifolds.
\end{enumerate}

Here is a brief outline to the paper. In Section \ref{sect-2}, we review briefly the representation theory we shall need
concerning $\mathcal{U}_\pm^\star$ submodules of $\otimes^kV^*$ and obtain an upper bound on the dimension of the space of
quadratic invariants for $\mathfrak{H}_\pm$ as a $\mathcal{U}_\pm^\star$ module. In Section
\ref{sect-3}, we turn to the geometric setting and study
$\nabla\Omega_\pm$. In Section \ref{sect-4}, we examine matters in the algebraic context and define projectors on the spaces
$W_{1,\pm}$, $W_{2,\pm}$, $U_{3,\pm}$, and $W_{4,\pm}$. In Section \ref{sect-5}, we fix the metric and vary the almost
(para)-complex structure to prove Theorem \ref{thm-1.1} and
Theorem
\ref{thm-1.4}. In Section \ref{sect-6}, we assume the (para)-complex structure to be integrable and vary the
metric to prove Theorem \ref{thm-1.2}. In Section
\ref{sect-7}, we use results of \cite{GH80} to establish Theorem \ref{thm-1.6}.

\section{Representation theory}\label{sect-2}
Let $(V,\langle\cdot,\cdot\rangle,J_\pm)$ be a para/pseudo-Hermitian space. Extend $\langle\cdot,\cdot\rangle$
to
$\otimes^kV$ so
\begin{equation}\label{eqn-2.a}
\langle(v_1\otimes\dots\otimes v_k),(w_1\otimes\dots\otimes w_k)\rangle:=\prod_{i=1}^k\langle
v_i,w_i\rangle\,.
\end{equation}
Equation~(\ref{eqn-2.a}) defines a non-degenerate symmetric bilinear form on $\otimes^kV$. We use $\langle\cdot,\cdot\rangle$
to identify
$V$ with $V^*$ and $\otimes^kV$ with $\otimes^kV^*$. If $\theta\in\otimes^kV^*$ and if $u\in\mathcal{U}_\pm^\star$,
the pull-back $u^*\theta\in\otimes^kV^*$ is defined by $u^*\theta(v_1,\dots,v_k):=\theta(uv_1,\dots,uv_k)$. Pull-back
defines a natural action of  $\mathcal{U}_\pm^\star$ on
$\otimes^kV^*$ which preserves the canonical inner product of Equation~(\ref{eqn-2.a}). Let $\xi$ be a
$\mathcal{U}_\pm^\star$-invariant subspace of
$\otimes^kV^*$; the natural action of
$\mathcal{U}_\pm^\star$ on
$\otimes^kV^*$ by pull-back makes $\xi$ into a $\mathcal{U}_\pm^\star$ submodule of $\otimes^kV^*$. One has:
\begin{lemma}\label{lem-2.1}
Let $(V,\langle\cdot,\cdot\rangle,J_\pm^0)$ be a para/pseudo-Hermitian vector space. Let $\xi$ be a
$\mathcal{U}_\pm^\star$ submodule of
$\otimes^kV^\ast$.
\begin{enumerate}\item $\langle\cdot,\cdot\rangle$ is non-degenerate on $\xi$.
\item There is an orthogonal direct sum decomposition $\xi=\eta_1\oplus\dots\oplus\eta_k$ where the $\eta_i$
are irreducible $\mathcal{U}_\pm^\star$-modules.
\item If $\xi_1$ and $\xi_2$ are inequivalent irreducible $\mathcal{U}_\pm^\star$ submodules of $\xi$, then $\xi_1\perp\xi_2$.
\item The multiplicity with which an irreducible representation appears in $\xi$ is independent of the decomposition in
{\rm(2)}.
\item If $\xi_1$ appears with multiplicity $1$ in $\xi$ and if $\eta$ is any $\mathcal{U}_\pm^\star$ submodule of $\xi$, then
either
$\xi_1\subset\eta$ or else $\xi_1\perp\eta$.
\item If $0\rightarrow\xi_1\rightarrow\xi\rightarrow\xi_2\rightarrow0$ is a short exact sequence of
$\mathcal{U}_\pm^\star$-modules, then
$\xi\approx\xi_1\oplus\xi_2$ as a $\mathcal{U}_\pm^\star$-module.
\end{enumerate}
\end{lemma}

\begin{proof} We shall establish Assertion (1) as this is the crucial property; the remaining assertions follow from Assertion
(1) using essentially the same arguments as those used in the positive definite setting; we refer to
\cite{BGN11} for a detailed exposition. For example, it is Assertion (1) which lets us define orthogonal projection; if $\xi$ is invariant
under the action of
$\mathcal{U}_\pm^\star$, then $\xi\cap\xi^\perp$ is a totally isotropic invariant subspace of $\otimes^kV^*$ and hence
$\xi\cap\xi^\perp=\{0\}$. Thus $\otimes^kV^*=\xi\oplus\xi^\perp$ and orthogonal projection on $\xi$ is given by the first
factor in this decomposition.

Suppose first $(V,\langle\cdot,\cdot\rangle,J_-^0)$ is a pseudo-Hermitian vector space of signature $(p,q)$. We prove
Assertion (1) for the smaller group $\mathcal{U}_-$; it then follows automatically for the larger group $\mathcal{U}_-^\star$.
Use the Gramm-Schmidt process to choose an orthogonal decomposition $V=V_+\oplus V_-$ which is $J_-^0$ invariant so $V_+$ is
spacelike and
$V_-$ is timelike. Let $T=\pm\id$ on $V_\pm$; $T\in\mathcal{U}_-$ since the decomposition is $J_-^0$ invariant. Let
$\{e_1,...,e_p\}$ be an orthonormal basis for
$V_-$ and let $\{e_{p+1},...,e_m\}$ be an orthonormal basis for $V_+$. Let $\{e^1,...,e^m\}$ be the corresponding orthonormal
dual basis for
$V^*$. Then $T^*(e^i)=\langle e^i,e^i\rangle e^i=\pm e^i$. If
$I=(i_1,...,i_k)$ is a multi-index, set
$e^I:=e^{i_1}\otimes...\otimes e^{i_k}$. The collection $\{e^I\}$ is an orthonormal basis for $\otimes^kV^*$ with:
\begin{eqnarray*}
T^*e^I&=&T^*(e^{i_1})\otimes...\otimes T^*(e^{i_k})
=\langle e^{i_1},e^{i_1}\rangle e^{i_1}\otimes...\otimes\langle e^{i_k},e^{i_k}\rangle e^{i_k}\\
&=&\langle e^I,e^I\rangle e^I=\pm e^I\,.
\end{eqnarray*}
Thus if $T^*w=w$, then $w$ is a spacelike vector in $\otimes^kV^*$ while if $T^*w=-w$, then $w$ is a timelike vector in
$\otimes^kV^*$. Let $\xi$ be a non-trivial $\mathcal{U}_-$ invariant subspace of $\otimes^kV^*$. Since $T\in\mathcal{U}_-$, $T$
preserves $\xi$. Decompose $\xi=\xi_+\oplus\xi_-$ into the $\pm1$ eigenspaces of $T^*$. Since $\xi_+$ is spacelike and $\xi_-$
is timelike, the metric on $\xi$ is non-degenerate and Assertion (1) follows in this framework.

The argument is a bit different in the para-Hermitian setting. Let $(V,\langle\cdot,\cdot\rangle,J_+^0)$ be a para-Hermitian
vector space. Find an orthogonal direct sum decomposition $V=V_+\oplus V_-$ where $V_+$ is spacelike, where $V_-$ is timelike,
and where $J_+^0:V_\pm\rightarrow V_\mp$. As before, let $T=\pm\id$ on $V_\pm$; $T$ does not belong to $\mathcal{U}_+$ but it
does belong to $\mathcal{U}_+^\star$. The remainder of the argument now follows as in the complex case; it is necessary to
assume $\xi$ is invariant under $\mathcal{U}_+^\star$ and not simply under $\mathcal{U}_+$ -- this is the crucial difference.
\end{proof}

\begin{remark}
\rm
Lemma \ref{lem-2.1} fails for the group $\mathcal{U}_+$ and it is for this reason that the decomposition of
$\mathfrak{H}_+$ has more factors as a $\mathcal{U}_+$ module than as a $\mathcal{U}_+^\star$ module. Let
$(V,\langle\cdot,\cdot\rangle,J_+^0)$ be a para-Hermitian vector space. Decompose $V=W_+\oplus W_-$ into the $\pm1$
eigenspaces of $J_+^0$. Then $W_\pm$ are totally isotropic subspaces of $V$ which are invariant under $\mathcal{U}_+$.
\end{remark}

Let $\xi$ be a $\mathcal{U}_\pm^\star$ submodule of $\otimes^kV^*$. We say that a symmetric inner product $\theta\in S^2(\xi^*)$ is
a {\it quadratic invariant} if
$\theta(\gamma x,\gamma y)=\theta(x,y)$ for all $\gamma\in \mathcal{U}_\pm^\star$ and for all $x,y\in\xi$; let
$S_{\mathcal{U}_\pm^\star}^2(\xi)$ be the space of all quadratic invariants. The following is well known -- see, for example, the
discussion in \cite{BGN11}. The proof follows exactly the same lines as in the positive definite setting given Lemma
\ref{lem-2.1} (1).
\begin{lemma}\label{lem-2.3}
Let $\xi$ be a $\mathcal{U}_\pm^\star$ submodule of
$\otimes^kV^*$. Suppose that $\xi_i$ are non-trivial $\mathcal{U}_\pm^\star$-modules so that $\xi_1\oplus\dots\oplus\xi_\ell$ is a
$\mathcal{U}_\pm^\star$ submodule of
$\xi$. Also suppose that
$\dim\{S_{\mathcal{U}_\pm^\star}^2(\xi)\}\le\ell$. Then:
\begin{enumerate}
\item $\xi=\xi_1\oplus\dots\oplus\xi_\ell$, $\xi_i\perp\xi_j$ for $i\ne j$, and $\dim\{S_{\mathcal{U}_\pm^\star}^2(\xi)\}=\ell$.
\item The modules $\xi_i$ are all irreducible and $\xi_i$ is not isomorphic to $\xi_j$ for $i\ne j$.
\end{enumerate}\end{lemma}

We now examine the space of quadratic invariants for the setting at hand.

\begin{lemma}\label{lem-2.4}
$\dim\{S^2_{\mathcal{U}_\pm^\star}(\mathfrak{H}_\pm)\}\le 4$.
\end{lemma}

\begin{proof} Since the original discussion in \cite{GH80} was in the positive definite setting, we shall provide full details.
Let $(V,\langle\cdot,\cdot\rangle,J_\pm^0)$ be a para/pseudo-Hermitian vector space and let $\xi$ be a $G$
submodule of $\otimes^kV^*$. A spanning set for the space of quadratic invariants if $G=\mathcal{O}$ or if $G=\mathcal{U}_-$ in
the positive definite setting is given in
\cite{W46} and in \cite{Fu58,Iw58}, respectively. The extension to the groups $\mathcal{U}_\pm^\star$ is straightforward
(see \cite{BGN11} for example). In brief, if $G=\mathcal{U}_\pm^\star$, everything is given by contraction of indices using the
inner product
$\langle\cdot,\cdot\rangle$ and the structure
$J_\pm^0$ where $J_\pm^0$ must appear an even number of times.  The following is a convenient formalism. We identify
$\theta$ with the corresponding quadratic function
$\theta(x):=\theta(x,x)$. We consider 3 distinct orthonormal bases
$\{e_{i_1}^1,e_{i_2}^2,e_{i_3}^3\}$ for $V$ which are indexed by
$\{i_1,i_2,i_3\}$, respectively, for $1\le i_1\le m$, $1\le i_2\le m$, and $1\le i_3\le m$. Let
$$\varepsilon_I=\langle e_{i_1}^1,e_{i_1}^1\rangle \langle e_{i_2}^2,e_{i_2}^2\rangle \langle
e_{i_3}^3,e_{i_3}^3\rangle=\pm1\,.
$$
We
consider a string $S$ of 6 symbols grouped into 2 monomials of 3 symbols where each index
$1$,
$2$,
$3$ appears twice and where some of the indices are decorated with $J_\pm^0$. Thus, for example, if
$S=(1,2;J_\pm^02)(1,3;J_\pm^03)$ and if $H_\pm\in\mathfrak{H}_\pm$, then the associated invariant $\mathcal{I}(S)$ is given by:
\begin{eqnarray*}
&&\mathcal{I}(S)(H_\pm):=\sum_{i_1=1}^m\sum_{i_2=1}^m\sum_{i_3=1}^m\varepsilon_IH_\pm(e_{i_1}^1,e_{i_2}^2;J_\pm^0e_{i_2}^2)
H_\pm(e_{i_1}^1,e_{i_3}^3;J_\pm^0e_{i_3}^3)\,.
\end{eqnarray*}
The space of quadratic invariants of $\mathfrak{H}_\pm$ is spanned by such invariants. We will stratify the invariants by the
number of times
$J_\pm^0$ appears; this gives rise to 2 basic cases each of which has 2 subcases.
\begin{enumerate}
\item General remarks.
\begin{enumerate}
\item We can
replace the basis $\{e_{i_1}^1\}$ by $\{J_\pm^0e_{i_1}^1\}$ and thereby replace $\varepsilon_I$ by $\mp\varepsilon_I$.
Thus
$\II((\dots,1,\dots,1,\dots))=\mp\II((\dots,J_\pm^01,\dots,J_\pm^01,\dots))$.
\item We need only consider strings where either a given index is undecorated or it is decorated
exactly once.
\item We may permute the bases. Thus\par
 $\II((1,2;3)(1,2;3))=\II((2,3;1)(2,3;1))$.
\item By Equation (\ref{eqn-1.a}),
$\II((\mu,\sigma;\star)(\star,\star;\star))=-\II((\sigma,\mu;\star)(\star,\star;\star))$\par$=
\pm\II((J_\pm^0\mu,J_\pm^0\sigma;\star)(\star,\star;\star))$.
\item By Equation  (\ref{eqn-1.b}),
$\II((\mu,J_\pm^0\sigma;\star)(\star,\star;\star))=\II((J_\pm^0\mu,\sigma;\star)(\star,\star;\star))$.
\end{enumerate}
\item $J_\pm^0$ does not appear. This gives rise to 3 invariants:
\begin{enumerate}
\item Each index appears in each variable:
\begin{enumerate}
\item $\psi_1:=\II((1,2;3)(1,2;3))$.
\item $\psi_2:=\II((1,2;3)(1,3;2))$.
\end{enumerate}
\item Only one index appears in both variables:
\begin{enumerate}
\item $\psi_3:=\II((1,2;1)(3,2;3))$.
\end{enumerate}\end{enumerate}
\item $J_\pm^0$ appears twice. This gives rise to another invariant:
\begin{enumerate}
\item Each index appears in each variable:
\begin{enumerate}
\item $\psi_4:=\II((1,J_\pm^02;J_\pm^03)(1,2;3))$.
\item $\II((1,J_\pm^02;3)(1,J_\pm^03;2))=\II((J_\pm^01,2;3)(J_\pm^01,2;3))$
\par\quad$=\mp\II((1,2;3)(1,2;3))=\mp\psi_1$.
\end{enumerate}
\item Only one index appears in both variables:
\begin{enumerate}
\item $\II((J_\pm^01,2;1)(J_\pm^03,2;3))=\II((1,J_\pm^02;1)(3,J_\pm^02;3))$
\par\quad
      $=\mp\II((1,2;1)(3,2;3))=\mp\psi_3$.
\end{enumerate}
\end{enumerate}
\end{enumerate}
We have enumerated all the possibilities and constructed 4 invariants.\end{proof}

\section{Geometric analysis}\label{sect-3}
If $(x^1,\dots,x^m)$ is a system of local coordinates on $M$, let
$\partial_{x_i}:=\frac{\partial}{\partial x_i}$.
\begin{lemma}\label{lem-3.1}
Let $(M,g,J_\pm)$ be an almost para/pseudo-Hermitian manifold. Then:
\begin{enumerate}
\item $\nabla\Omega_\pm\in\mathfrak{H}_\pm$.
\item $\nabla\Omega_\pm(x,y;z)=g(x,(\nabla_zJ_\pm)y)
=g(x,\nabla_zJ_\pm y)-g(x,J_\pm\nabla_zy)$\par
\noindent$\phantom{\nabla\Omega_\pm(x,y;z)}=g(x,\nabla_zJ_\pm y)+g(J_\pm x,\nabla_zy)$.
\end{enumerate}
\end{lemma}

\begin{proof}
Since $\Omega_\pm\in C^\infty(\Lambda^2)$,
$\nabla\Omega_\pm\in C^\infty(\Lambda^2\otimes
V^*)$. We prove Assertion (1) by studying the action of $J_\pm^*$:
\medbreak\quad
$\nabla\Omega_\pm(J_\pm x,J_\pm y;z)$
\smallbreak\qquad
$=zg(J_\pm x,J_\pm J_\pm y)-g(\nabla_zJ_\pm x,J_\pm J_\pm y)-g(J_\pm x,J_\pm \nabla_zJ_\pm y)$
\smallbreak\qquad
$=\mp zg(x,J_\pm y)\mp g(\nabla_zJ_\pm x,y)\pm g(x,\nabla_zJ_\pm y)$
\smallbreak\qquad
$=\mp zg(x,J_\pm y)\mp zg(J_\pm x,y)\pm g(J_\pm x,\nabla_zy)$
$\pm zg(x,J_\pm y)\mp g(\nabla_zx,J_\pm y)$
\smallbreak\qquad
$=\pm zg(x,J_\pm y)\mp g(x,J_\pm \nabla_zy)\mp g(\nabla_zx,J_\pm y)$
$=\pm\nabla\Omega_\pm(x,y;z)$.
\medbreak\noindent We use the fact that $\nabla g=0$ to prove Assertion (2) by computing:
\medbreak\quad
$\nabla_z\Omega_\pm(x,y)=zg(x,J_\pm y)-g(\nabla_zx,J_\pm y)-g(x,J_\pm\nabla_zy)$
\smallbreak\quad
$=zg(x,J_\pm y)-g(\nabla_zx,J_\pm y)-g(x,\nabla_zJ_\pm y)+g(x,\nabla_zJ_\pm y)-g(x,J_\pm\nabla_zy)$
\smallbreak\quad
$=(\nabla_zg)(x,J_\pm y)+g(x,\nabla_zJ_\pm y)-g(x,J_\pm\nabla_zy)$
\smallbreak\quad
$=g(x,\nabla_zJ_\pm y)-g(x,J_\pm\nabla_zy)=g(x,\nabla_zJ_\pm y)+g(J_\pm x,\nabla_zy)$.
\end{proof}

Let $g(x,y;z):=zg(x,y)$. We continue our study and assume $J_\pm$ is integrable:

\begin{lemma}\label{lem-3.2}
Let $(M,g,J_\pm)$ be a para/pseudo-Hermitian manifold. Then:
 \begin{enumerate}
\item $\nabla\Omega_\pm(\partial_{x_i},\partial_{x_j};\partial_{x_k})=\textstyle\frac12\{
g(\partial_{x_i},\partial_{x_k};J_\pm \partial_{x_j})-g(\partial_{x_j},\partial_{x_k};J_\pm \partial_{x_i})$
\par$+g(J_\pm\partial_{x_i},\partial_{x_k};\partial_{x_j})-g(J_\pm \partial_{x_j},\partial_{x_k};\partial_{x_i})\}$.
\smallbreak\item $\nabla\Omega_\pm(M,g,J_\pm)\in U_{3,\pm}$.
\smallbreak\item $\nabla\Omega_\pm(M,e^{2f}g,J_\pm)=e^{2f}\{\nabla\Omega_\pm(M,g,J_\pm)-\sigma_{\pm,g}(df)\}$.
\smallbreak\item $W_{4,\pm}\subset U_{3,\pm}$.
\end{enumerate}
\end{lemma}

\begin{proof}
Since $J_\pm$ is integrable, we may choose coordinates so $J_\pm\partial_{x_i}\in\{\partial_{x_1},...,\partial_{x_m}\}$.
Let $x=\partial_{x_i}$, $y=\partial_{x_j}$, and $z=\partial_{x_k}$. We may apply Lemma \ref{lem-3.1} and the Koszul
formula for the Christoffel symbols in a coordinate frame to see:
\begin{eqnarray*}
\nabla_z\Omega_\pm(x,y)&=&g(x,\nabla_zJ_\pm y)+g(J_\pm x,\nabla_zy)\\
&=&\textstyle\frac12\{g(x,z;J_\pm y)+g(x,J_\pm y;z)-g(z,J_\pm y;x)\}\\
&+&\textstyle\frac12\{g(J_\pm x,z;y)+g(J_\pm x,y;z)-g(z,y;J_\pm x)\}\,.
\end{eqnarray*}
Assertion (1) now follows from the identity:
$$g(x,J_\pm y;z)+g(J_\pm x,y;z)=z\{g(x,J_\pm y)+g(J_\pm x,y)\}=0\,.$$

We prove Assertion (2) by checking that $\nabla\Omega_\pm$
satisfies the defining relation for $U_{3,\pm}$ in this instance. We use Assertion (1) to compute:
\begin{eqnarray*}
&&\nabla\Omega_\pm(x,J_\pm y;J_\pm z)\\
&=&\textstyle\frac12\{g(x,J_\pm  z;J_\pm J_\pm  y)-g(J_\pm  y,J_\pm  z;J_\pm x)\}\\
&+&\textstyle\frac12\{g(J_\pm x,J_\pm
z;J_\pm  y)-g(J_\pm J_\pm  y,J_\pm  z;x)\}\\
&=&\textstyle\frac12\{\pm g(x,J_\pm z;y)\pm g(y,z;J_\pm x)\mp g(x,z;J_\pm  y)\pm g(J_\pm y,z;x)\}\\
&=&\textstyle\frac12\{\mp g(J_\pm x,z;y)\pm g(y,z;J_\pm x)\mp g(x,z;J_\pm  y)\pm g(J_\pm y,z;x)\}\\
&=&\mp\nabla_\pm(x,y;z)\,.
\end{eqnarray*}

We also use Assertion (1) to prove Assertion (3) by checking:
\begin{eqnarray*}
&&\nabla\Omega_{\pm,e^{2f}g}(x,y;z)\\
&=&e^{2f}\{\nabla\Omega_{\pm,g}(x,y;z)+df(J_\pm y)g(x,z)-df(J_\pm x)g(y,z)\\
&+&df(y)g(J_\pm x,z)-df(x)g(J_\pm y,z)\}\\
&=&e^{2f}\{\nabla\Omega_{\pm,g}(x,y;z)-\sigma_{\pm,g}(df)(x,y;z)\}\,.
\end{eqnarray*}

Let $(V,\langle\cdot,\cdot\rangle,J_\pm^0)$ be a para/pseudo-Hermitian vector space. Let $f$ be a smooth
function on $V$ and consider the manifold $(M,g,J_\pm):=(V,e^{2f}\langle\cdot,\cdot\rangle,J_\pm^0)$. We apply Assertion (2)
and Assertion (3) to prove Assertion (4) by checking:
\medbreak\noindent\hfill
$e^{2f}\sigma_{\pm,\langle\cdot,\cdot\rangle}(df)=-\nabla\Omega_{\pm,e^{2f}\langle\cdot,\cdot\rangle}\in
U_{3,\pm}$.\hfill
\end{proof}

\section{Algebraic considerations}\label{sect-4}
We now turn our attention to purely algebraic considerations. For the remainder of this section, let
$(V,\langle\cdot,\cdot\rangle,J_\pm^0)$ be a para/pseudo-Hermitian vector space.

\begin{definition}
\rm
Let $H_\pm\in\mathfrak{H}_\pm$.\begin{enumerate}
\item
$(\pi_{1,\pm}H_\pm)(x,y;z):=\frac16\big\{H_\pm(x,y;z)+H_\pm(y,z;x)+H_\pm(z,x;y)$
\smallbreak\qquad
$\pm H_\pm(x,J_\pm^0 y;J_\pm^0 z)\pm H_\pm(y,J_\pm^0z;J_\pm^0x)\pm H_\pm(z,J_\pm^0x;J_\pm^0y)\big\}$.
\smallbreak\item
$(\pi_{2,\pm}H_\pm)(x,y;z):=\frac16\big\{2H_\pm(x,y;z)-H_\pm(y,z;x)-H_\pm(z,x;y)$
\smallbreak\qquad
$\pm 2H_\pm(x,J_\pm^0 y;J_\pm^0 z)\mp H_\pm(y,J_\pm^0z;J_\pm^0x)\mp H_\pm(z,J_\pm^0x;J_\pm^0y)\big\}$.
\smallbreak\item $(\pi_{3,\pm}H_\pm)(x,y;z):=\textstyle\frac12\{H_\pm(x,y;z)\mp H_\pm(x,J_\pm^0 y;J_\pm^0 z)\}$.
\smallbreak\item $\pi_{4,\pm}:=\pm\textstyle\frac1{m-2}\sigma_\pm (J_\pm^0)^*\tau_1$.
\end{enumerate}\end{definition}

\begin{lemma}\label{lem-4.2}
\ \begin{enumerate}
\item $\pi_{1,\pm}$ is a projection from $\mathfrak{H}_\pm$ onto $W_{1,\pm}$.
\smallbreak\item $\pi_{2,\pm}$ is a projection from $\mathfrak{H}_\pm$ onto $W_{2,\pm}$.
\smallbreak\item $\pi_{3,\pm}$ is a projection from $\mathfrak{H}_\pm $ onto $U_{3,\pm}$.
\smallbreak\item $\pi_{4,\pm}$ is a projection from $\mathfrak{H}_\pm$ onto $W_{4,\pm}$.
\end{enumerate}\end{lemma}

\begin{proof} Set:
\begin{eqnarray*}
&&(\kappa_\pm^1H_\pm)(x,y;z):=\pm H_\pm(x,J_\pm^0y;J_\pm^0 z),\\
&&(\kappa_\pm^2H_\pm)(x,y;z):=H_\pm(y,z;x)+H_\pm(z,x;y)\\
&&\qquad\qquad\qquad\pm H_\pm(y,J_\pm^0z;J_\pm^0x)\pm H_\pm(z,J_\pm^0x;J_\pm^0y)\,.
\end{eqnarray*}
We may use Equation~(\ref{eqn-1.a}) and Equation~(\ref{eqn-1.b}) to see that $\kappa_\pm^1H_\pm$, and $\kappa_\pm^2H_\pm$ are
anti-symmetric in
the first two arguments. We show that $\kappa_\pm^1\mathfrak{H}_\pm\subset\mathfrak{H}_\pm$ and that
$\kappa_\pm^2\mathfrak{H}_\pm\subset\mathfrak{H}_\pm$ by checking:
\medbreak\quad
$(\kappa_\pm^1H_\pm)(J_\pm^0x,J_\pm^0y;z)=H_\pm(J_\pm^0x,J_\pm^0J_\pm^0y;J_\pm^0z)=\pm H_\pm(x,J_\pm^0y;J_\pm^0z)$
\smallbreak\qquad\quad
$=\pm\kappa_\pm^1H_\pm(x,y;z)$,
\medbreak\quad
$(\kappa_\pm^2H_\pm)(J_\pm^0x,J_\pm^0y;z)=H_\pm(J_\pm^0y,z;J_\pm^0x)+H_\pm(z,J_\pm^0x;J_\pm^0y)$
\smallbreak\qquad\quad
$\pm H_\pm(J_\pm^0y,J_\pm^0z;J_\pm^0J_\pm^0x)\pm H_\pm(z,J_\pm^0J_\pm^0x;J_\pm^0J_\pm^0y)$
\smallbreak\qquad
$=H_\pm(y,J_\pm^0z;J_\pm^0x)+H_\pm(z,J_\pm^0x;J_\pm^0y)
\pm H_\pm(y,z;x)\pm H_\pm(z,x;y)$
\smallbreak\qquad
$=\pm(\kappa_\pm^2H_\pm)(x,y;z)$.
\medbreak\noindent We see $\pi_{1,\pm}\mathfrak{H}_\pm\subset\mathfrak{H}_\pm$,
$\pi_{2,\pm}\mathfrak{H}_\pm\subset\mathfrak{H}_\pm$, and
$\pi_{3,\pm}\mathfrak{H}_\pm\subset\mathfrak{H}_\pm$ by expressing:
\begin{eqnarray*}
&&\pi_{1,\pm}=\textstyle\frac16\{\id+\kappa_\pm^1+\kappa_\pm^2\},\quad
  \pi_{2,\pm}=\textstyle\frac16\{2(\id+\kappa_\pm^1)-\kappa_\pm^2\},\\
&&\pi_{3,\pm}=\textstyle\frac12\{\id-\kappa_\pm^1\}\,.
\end{eqnarray*}

Let $H_\pm\in\mathfrak{H}_\pm$. We verify $\pi_{1,\pm}H_\pm\in W_{1,\pm}$, that
$\pi_{2,\pm}H_\pm\in W_{2,\pm}$, and that $\pi_{3,\pm}H_\pm\in U_{3,\pm}$ by checking that the
defining relations are satisfied in each case:
\medbreak\quad
$(\pi_{1,\pm}H_\pm)(x,z;y):=
\frac16\big\{H_\pm(x,z;y)+H_\pm(z,y;x)+H_\pm(y,x;z)$
\smallbreak\qquad\quad
$\pm H_\pm(x,J_\pm^0 z;J_\pm^0 y)\pm H_\pm(z,J_\pm^0y;J_\pm^0x)\pm H_\pm(y,J_\pm^0x;J_\pm^0z)\big\}$
\smallbreak\qquad
$=-\pi_{1,\pm}H_\pm(x,y;z)$,
\medbreak\quad
$(\pi_{2,\pm}H_\pm)(x,y;z)+(\pi_{2,\pm}H_\pm)(y,z;x)+(\pi_{2,\pm}H_\pm)(z,x;y)$
\smallbreak\qquad
$=\frac16\big\{2H_\pm(x,y;z)-H_\pm(y,z;x)-H_\pm(z,x;y)$
\smallbreak\qquad\quad
$\pm 2H_\pm(x,J_\pm^0 y;J_\pm^0 z)\mp H_\pm(y,J_\pm^0z;J_\pm^0x)
\mp H_\pm(z,J_\pm^0x;J_\pm^0y)$
\smallbreak\qquad\quad
$+2H_\pm(y,z;x)-H_\pm(z,x;y)-H_\pm(x,y;z)$
\smallbreak\qquad\quad
$\pm 2H_\pm(y,J_\pm^0 z;J_\pm^0 x)\mp H_\pm(z,J_\pm^0x;J_\pm^0y)
\mp H_\pm(x,J_\pm^0y;J_\pm^0z)$
\smallbreak\qquad\quad
$+2H_\pm(z,x;y)-H_\pm(x,y;z)-H_\pm(y,z;x)$
\smallbreak\qquad\quad
$\pm 2H_\pm(z,J_\pm^0 x;J_\pm^0 y)\mp H_\pm(x,J_\pm^0y;J_\pm^0z)
\mp H_\pm(y,J_\pm^0z;J_\pm^0x)\big\}$
$=0$,
\medbreak\quad
$(\pi_{3,\pm}H_\pm)(x,J_\pm^0y;J_\pm^0z)$
\smallbreak\qquad
$=\frac12\{H_{\pm}(x,J_\pm^0y;J_\pm^0z)\mp H_\pm (x,J_\pm^0J_\pm^0y;J_\pm^0J_\pm^0z)\}$
\smallbreak\qquad
$=\frac12\{H_{\pm}(x,J_\pm^0y;J_\pm^0z)\mp H_\pm(x,y;z)\}$
\smallbreak\qquad
$=\mp\frac12\{H(x,y;z)\mp H(x,J_\pm^0y;J_\pm^0z)\}$
$=\mp(\pi_{3,\pm}H_\pm)(x,y;z)$.

\medbreak Let $H_{1,\pm}\in W_{1,\pm}$, let $H_{2,\pm}\in W_{2,\pm}$, and let $H_{3,\pm}\in
U_{3,\pm}$. We complete the proof of Assertion (1), of Assertion (2), and of Assertion (3) by verifying:
\medbreak\quad
$(\pi_{1,\pm}H_{1,\pm})(x,y;z)=\frac16\big\{H_{1,\pm}(x,y;z)+H_{1,\pm}(y,z;x)+H_{1,\pm}(z,x;y)$
\smallbreak\qquad\quad
$\pm H_{1,\pm}(x,J_\pm^0 y;J_\pm^0 z)\pm H_{1,\pm}(y,J_\pm^0z;J_\pm^0x)\pm H_{1,\pm}(z,J_\pm^0x;J_\pm^0y)\big\}$
\smallbreak\qquad
$=\frac16\big\{H_{1,\pm}(x,y;z)-H_{1,\pm}(y,x;z)-H_{1,\pm}(x,z;y)$
\smallbreak\qquad\quad
$\mp H_{1,\pm}(J_\pm^0 z,J_\pm^0 y;x)\mp H_{1,\pm}(J_\pm^0x,J_\pm^0z;y)\mp H_{1,\pm}(J_\pm^0y,J_\pm^0x;z)\big\}$
\smallbreak\qquad
$=\frac16\big\{H_{1,\pm}(x,y;z)+H_{1,\pm}(x,y;z)+H_{1,\pm}(x,y;z)$
\smallbreak\qquad\quad
$- H_{1,\pm}(z,y;x)-H_{1,\pm}(x,z;y)-H_{1,\pm}(y,x;z)\big\}$
$=H_{1,\pm}(x,y;z)$,
\medbreak\quad
$(\pi_{2,\pm}H_{2,\pm})(x,y;z)=\frac16\big\{2H_{2,\pm}(x,y;z)-H_{2,\pm}(y,z;x)-H_{2,\pm}(z,x;y)$
\smallbreak\qquad\quad
$\pm2H_{2,\pm}(x,J_\pm^0 y;J_\pm^0 z)\mp H_{2,\pm}(y,J_\pm^0z;J_\pm^0x)
\mp H_{2,\pm}(z,J_\pm^0x;J_\pm^0y)\big\}$
\smallbreak\qquad
$=\frac16\big\{3H_{2,\pm}(x,y;z)\pm H_{2,\pm}(J_\pm^0x,y;J_\pm^0z)\pm H_{2,\pm}(x,J_\pm^0y;J_\pm^0z)$
\smallbreak\qquad\quad
$\mp H_{2,\pm}(y,J_\pm^0z;J_\pm^0x)
\mp H_{2,\pm}(J_\pm^0z,x;J_\pm^0y)\big\}$
\smallbreak\qquad
$=\frac16\big\{3H_{2,\pm}(x,y;z)\mp H_{2,\pm}(J_\pm^0z,J_\pm^0x;y)\mp H_{2,\pm}(y,J_\pm^0z;J_\pm^0x)$
\smallbreak\qquad\quad
$\mp H_{2,\pm}(J_\pm^0y,J_\pm^0z;x)\mp H_{2,\pm}(J_\pm^0z,x;J_\pm^0y)$
\smallbreak\qquad\quad
$\mp H_{2,\pm}(y,J_\pm^0z;J_\pm^0x)
\mp H_{2,\pm}(J_\pm^0z,x;J_\pm^0y)\big\}$
\smallbreak\qquad
$=\frac16\{3H_{2,\pm}(x,y;z)-H_{2,\pm}(z,x;y)-H_{2,\pm}(y,z;x)$
\smallbreak\qquad\quad
$\mp 2H_{2,\pm}(J_\pm^0y,z;J_\pm^0x)\mp2H_{2,\pm}(z;J_\pm^0x;J_\pm^0y)\}$
\smallbreak\qquad
$=\frac16\{4H_{2,\pm}(x,y;z)\pm 2H_{2,\pm}(J_0^\pm x,J_0^\pm y;z)\}=H_{2,\pm}(x,y;z)$,
\medbreak\quad
$(\pi_{3,\pm}H_{3,\pm})(x,y;z)=\frac12\{H_{3,\pm}(x,y;z)\mp H_{3,\pm}(x,J_\pm^0y;J_\pm^0z\}$
\smallbreak\qquad
$=\frac12\{H_{3,\pm}(x,y;z)+H_{3,\pm}(x,y;z)\}=H_{3,\pm}(x,y;z)$.

\medbreak We now turn to the final assertion.
We
compute:
\begin{eqnarray*}
&&\tau_1(\sigma_\pm(\phi))(x)=
\varepsilon^{ij}\{\phi(J_\pm^0x)\langle e_i,e_j\rangle-\phi(J_\pm^0e_i)\langle x,e_j\rangle\}\\
&&\qquad+\varepsilon^{ij}\{\phi(x)\langle J_\pm^0e_i,e_j\rangle-\phi(e_i)\langle J_\pm^0x,e_j\rangle\}\\
&=&m\phi(J_\pm^0x)-\phi(J_\pm^0x)+\trace(J_\pm^0)\phi(x)-\phi(J_\pm^0x)\\
&=&(m-2)((J_\pm^0)^*\phi)(x)+\trace(J_\pm^0)\phi(x)\,.
\end{eqnarray*}
Since $\trace(J_\pm^0)=0$, we have
$\tau_1\sigma_\pm=(m-2)(J_\pm^0)^*$. It is immediate that $\pi_{4,\pm}$ takes values in $W_{4,\pm}$. We
complete the proof by checking:
\medbreak\qquad
$\pi_{4,\pm}\sigma_\pm\phi=
  \pm\frac1{m-2}(\sigma_\pm (J_\pm^0)^*\tau_1)(\sigma_\pm\phi)$
$=\pm\sigma_\pm(J_\pm^0)^*(J_\pm^0)^*\phi=\sigma_\pm\phi$.
\end{proof}

We examine these modules further:
\begin{lemma}\label{lem-4.3}
\ \begin{enumerate}
\item $W_{1,\pm}+W_{2,\pm}\subset\ker\pi_{3,\pm}$.
\smallbreak\item $W_{1,\pm}\cap W_{2,\pm}=\{0\}$.
\smallbreak\item $W_{1,\pm}\oplus W_{2,\pm}\oplus W_{3,\pm}\oplus W_{4,\pm}$
is a $\mathcal{U}_\pm^\star$ submodule of $\mathfrak{H}_\pm$.
\end{enumerate}\end{lemma}

\begin{proof} Suppose first that $H_{1,\pm}\in W_{1,\pm}$. Then
\medbreak\quad
$\pi_{3,\pm}H_{1,\pm}(x,y;z)=\frac12\{H_{1,\pm}(x,y;z)\mp H_{1,\pm}(x,J_\pm^0y;J_\pm^0z)\}$
\smallbreak\qquad
$=\frac12\{H_{1,\pm}(x,y;z)\mp H_{1,\pm}(J_\pm^0x,y;J_\pm^0z)\}$
\smallbreak\qquad
$=\frac12\{H_{1,\pm}(x,y;z)\pm H_{1,\pm}(J_\pm^0x,J_\pm^0z;y)\}$
\smallbreak\qquad
$=\frac12\{H_{1,\pm}(x,y;z)+H_{1,\pm}(x,z;y)\}=0$.
\medbreak\noindent Next suppose that $H_{2,\pm}\in W_{2,\pm}$. We have
\medbreak\quad
$\pi_{3,\pm}H_{2,\pm}(x,y;z)=\frac12\{H_{2,\pm}(x,y;z)\mp H_{2,\pm}(x,J_\pm^0y;J_\pm^0z)\}$
\smallbreak\qquad
$=\frac12\{H_{2,\pm}(x,y;z)\pm H_{2,\pm}(J_\pm^0y,J_\pm^0z;x)\pm H_{2,\pm}(J_\pm^0z,x;J_\pm^0y)\}$
\smallbreak\qquad
$=\frac12\{H_{2,\pm}(x,y;z)+H_{2,\pm}(y,z;x)\pm H_{2,\pm}(z,J_\pm^0x;J_\pm^0y)\}$
\smallbreak\qquad
$=\frac12\{H_{2,\pm}(x,y;z)+H_{2,\pm}(y,z;x)\mp H_{2,\pm}(J_\pm^0x,J_\pm^0y;z)$
\smallbreak\qquad\qquad
$\mp H_{2,\pm}(J_\pm^0y,z;J_\pm^0x)\}$
\smallbreak\qquad
$=\frac12\{H_{2,\pm}(y,z;x)\mp H_{2,\pm}(J_\pm^0y,z;J_\pm^0x)\}$
\smallbreak\qquad
$=-\frac12\{H_{2,\pm}(z,y;x)\mp H_{2,\pm}(z,J_\pm^0y;J_\pm^0x)\}$
$=-\pi_{3,\pm}H_{2,\pm}(z,y;x)$.
\medbreak\noindent
This shows that
\begin{eqnarray*}
\pi_{3,\pm}H_{2,\pm}(x,y;z)&=&-\pi_{3,\pm}H_{2,\pm}(y,x;z)=\pi_{3,\pm}H_{2,\pm}(z,x;y)\\
&=&-\pi_{3,\pm}H_{2,\pm}(x,z;y)\,.
\end{eqnarray*}
Consequently $H_{1,\pm}:=\pi_{3,\pm}H_{2,\pm}\in W_{1,\pm}$. Thus:
$$
\pi_{3,\pm}H_{2,\pm}=
  \pi_{3,\pm}\pi_{3,\pm}H_{2,\pm}=\pi_{3,\pm}H_{1,\pm}=0\,.$$

Let $H_\pm\in W_{1,\pm}\cap W_{2,\pm}$. We establish Assertion (2) by checking:
\begin{eqnarray*}
0&=&H_\pm(x,y;z)+H_\pm(y,z;x)+H_\pm(z,x;y)\\
&=&H_\pm(x,y;z)-H_\pm(y,x;z)-H_\pm(x,z;y)\\
&=&3H_\pm(x,y;z)\,.
\end{eqnarray*}

If $\pi_\pm:\mathfrak{H}_\pm\rightarrow\mathfrak{H}_\pm$ satisfies $\pi_\pm^2=\pi_\pm$, then
Lemma \ref{lem-2.1} shows
$$\mathfrak{H}_\pm=\ker(\pi_\pm)\oplus\operatorname{Range}(\pi_\pm)\,.$$
By Lemma \ref{lem-4.2}, we can apply this observation to $\pi_{3,\pm}$ and to $\pi_{4,\pm}$. By Assertion (1) and by
Assertion (2),
$$
  W_{1,\pm}\cap W_{2,\pm}=\{0\}\quad\text{so}\quad
  W_{1,\pm}\oplus W_{2,\pm}\subset\ker(\pi_{3,\pm})\,.
$$
By Lemma \ref{lem-4.2},
we have
$U_{3,\pm}=\operatorname{Range}(\pi_{3,\pm})$. Consequently
$$W_{1,\pm}\oplus W_{2,\pm}\oplus U_{3,\pm}$$
is a submodule of $\mathfrak{H}_\pm$. By Lemma \ref{lem-4.2},
$$W_{4,\pm}=\operatorname{Range}(\pi_{4,\pm})\subset U_{3,\pm}.$$
Since
$W_{4,\pm}=\pi_{4,\pm}U_{3,\pm}$, $W_{3,\pm}\oplus W_{4,\pm}$
is a $\mathcal{U}_\pm^\star$ submodule of $U_{3,\pm}$.
\end{proof}

\section{Varying the almost (para)-complex structure}\label{sect-5}
Fix a background almost para/pseudo-Hermitian manifold
$(M,g,J_\pm)$ and a point $P$ of $M$ for the remainder of Section \ref{sect-5}. Let $\mathcal{O}(M)$ be the fiber bundle whose
fibre over a point $Q$ of $M$ is the associated structure group
$\mathcal{O}(T_QM,g_Q)$. The Lie algebra $\mathfrak{o}$ of $\mathcal{O}$
is the vector space of all matrices which are skew-adjoint with respect to the inner product. Let
$\vartheta\in\mathfrak{o}_P\otimes T_P^*$M. Let $\varTheta$ be a smooth section to $\mathcal{O}(M)$ so that
$\varTheta(P)=\id$, so that $\varTheta=\id$ off a neighborhood of $P$, and so that
$d\varTheta=\vartheta$. Let:
$$
J_\pm^{\varTheta}:=\varTheta^{-1}J_\pm\varTheta\,.
$$
Since $\varTheta$ takes values in $\mathcal{O}$, $(M,g,J_\pm^\varTheta)$ is an almost para/pseudo-Hermitian manifold as
well. Define:
$$\Xi_\pm(\vartheta)(x,y;z):=g(x,(-\vartheta(z)J_\pm+J_\pm\vartheta(z))y)(P)\,.$$

\goodbreak\begin{lemma}\label{lem-5.1}
Adopt the notation established above.
\begin{enumerate}
\item $\left\{\nabla\Omega_\pm(M,g,J_\pm^\varTheta)(x,y;z)-\nabla\Omega_\pm(M,g,J_\pm)(x,y;z)\right\}(P)
$\par$=\Xi_\pm(d\vartheta)(x,y;z)$.
\smallbreak\item $\Xi_\pm$ is a $\mathcal{U}_\pm^\star$ module morphism from $\mathfrak{o}\otimes
V^*\otimes\chi$ to $\mathfrak{H}_\pm$.
\smallbreak\item If $m\ge6$, then $\pi_{1,\pm}\{\Xi_\pm(\mathfrak{o})\}\ne\{0\}$ and
$\pi_{3,\pm}\{\Xi_\pm(\mathfrak{o})\}\cap W_{3,\pm}\ne\{0\}$.
\smallbreak\item $\pi_{2,\pm}\{\Xi_\pm(\mathfrak{o})\}\ne\{0\}$ and
$\pi_{4,\pm}\{\Xi_\pm(\mathfrak{o})\}\ne\{0\}$.
\end{enumerate}
\end{lemma}

\begin{proof}Since $\varTheta(P)=\id$, $(J_\pm^\varTheta-J_\pm)(P)=0$. We use Lemma \ref{lem-3.1} to
prove Assertion (1) by computing:
\medbreak\qquad\qquad\qquad
$\Omega_\pm(M,g,J_\pm^\varTheta)(P)-\Omega_\pm(M,g,J_\pm)(P)$
\smallbreak\quad\qquad\qquad
$=g(x,\{\nabla_z(J_\pm^\varTheta-J_\pm)-(J_\pm^\varTheta-J_\pm)\nabla_z\}y)(P)$
\smallbreak\quad\qquad\qquad
$=g(x,\{z(J_\pm^\varTheta-J_\pm)\}y)(P)=g(x,\{z(\varTheta^{-1}J_\pm\varTheta-J_\pm)\}y)(P)$
\smallbreak\quad\qquad\qquad
$=g(x,\{-z(\varTheta)J_\pm+J_\pm z(\varTheta)\}y)(P)$.
\medbreak

Assertion (2) is an immediate consequence of Assertion (1). The proof of Assertions (3) and (4) is a purely algebraic
computation. Introduce an orthonormal basis $\{e_1,\dots,e_{\bar m},f_1,\dots,f_{\bar m}\}$ for $V$ so
$$J_\pm:e_i\rightarrow f_i\quad\text{and}\quad J_\pm:f_i\rightarrow\pm e_i\,.$$
We set $\varepsilon_i:=\langle e_i,e_i\rangle$. Define $\vartheta_0\in\mathfrak{o}$ by setting:
$$
\vartheta_0e_i=\left\{\begin{array}{rrr}
\varepsilon_2e_2&\text{if}&i=1\\
-\varepsilon_1e_1&\text{if}&i=2\\
0&\text{if}&i>2\end{array}\right\}\quad\text{and}\quad
\vartheta_0f_i=\left\{\begin{array}{rrr}
0&\text{if}&i=1\\0&\text{if}&i=2\\0&\text{if}&i>2\end{array}\right\}.
$$
Suppose first that $m\ge6$. We set $\vartheta=\vartheta_0\otimes e^3$. Choose $\alpha\in C^\infty(M)$ to be compactly supported
near $P$ with $d\alpha(P)=dx^3$. If
$\varepsilon_1=\varepsilon_2$, then the corresponding
$\varTheta$ may be taken to be:
$$\varTheta\partial_{x_i}=\left\{\begin{array}{rll}
\cos(\alpha)e_1+\varepsilon_2\sin(\alpha)e_2&\text{if}&i=1\\
-\varepsilon_2\sin(\alpha)e_1+\cos(\alpha)e_2&\text{if}&i=2\\
e_i&\text{if}&i\ge3\end{array}\right\}\quad\text{and}\quad
\varTheta\partial_{y_i}=\partial_{y_i}\ \forall\ i\,,$$
whereas if $\varepsilon_1=-\varepsilon_2$, then $\varTheta$ may be taken to be:
$$\varTheta\partial_{x_i}=\left\{\begin{array}{rll}
\cosh(\alpha)e_1+\varepsilon_2\sinh(\alpha)e_2&\text{if}&i=1\\
\varepsilon_2\sinh(\alpha)e_1+\cosh(\alpha)e_2&\text{if}&i=2\\
e_i&\text{if}&i\ge3\end{array}\right\}\quad\text{and}\quad
\varTheta\partial_{y_i}=\partial_{y_i}\ \forall\ i\,.$$
Set
$H_\pm:=\Xi_\pm(\vartheta_0\otimes e^3)$. The non-zero components of
$H_\pm$ are determined by:
$$
H_\pm(f_2,e_1;e_3)=\mp1\quad\text{and}\quad
H_\pm(f_1,e_2;e_3)=\pm1\,.
$$
Clearly $\tau_1H_\pm=0$; thus $\pi_{3,\pm}H_\pm\in W_{3,\pm}$. We prove Assertion (3) by computing:
$$
\pi_{1,\pm}H_\pm(f_2,e_1;e_3)=\mp\textstyle\frac16\quad\text{and}\quad
\pi_{3,\pm}H_\pm(f_2,e_1;e_3)=\mp\textstyle\frac12\,.
$$
Next we clear the previous notation and let $H_\pm=\Xi_\pm(\vartheta_0\otimes e^2)$; here we need to have $d\alpha(P)=dx^2$.
The non-zero components of
$H_\pm$ are determined by:
$$
H_\pm(f_2,e_1;e_2)=\mp1\quad\text{and}\quad
H_\pm(f_1,e_2;e_2)=\pm1\,.
$$
Since $\tau_1(H_\pm)=\pm\varepsilon_2$, the component of $H_\pm$ in $W_{4,\pm}$ is non-zero.
We complete the proof of Assertion (4) by checking:
\medbreak\qquad
$(\pi_{2,\pm}H_\pm)(f_2,f_1;f_2):=\frac16\big\{2H_\pm(f_2,f_1;f_2)-H_\pm(f_1,f_2;f_2)-H_\pm(f_2,f_2;f_1)$
\smallbreak\qquad\qquad
$\pm 2H_\pm(f_2,J_\pm^0 f_1;J_\pm^0 f_2)\mp H_\pm(f_1,J_\pm^0f_2;J_\pm^0f_2)\mp H_\pm(f_2,J_\pm^0f_2;J_\pm^0f_1)\big\}$
\smallbreak\qquad\quad
$=\frac16\{0-0-0-2-1+0\}=-\frac12$.
\end{proof}

\noindent{\bf Proof of Theorem \ref{thm-1.4}.} Let $m\ge6$. By Lemma \ref{lem-5.1}, $W_{i,\pm}$ are non-trivial modules for
$1\le i\le 4$. By Lemma \ref{lem-4.3}, $W_{1,\pm}\oplus W_{2,\pm}\oplus W_{3,\pm}\oplus W_{4,\pm}$
is a $\mathcal{U}_\pm^\star$ submodule of $\mathfrak{H}_\pm$. By Lemma \ref{lem-2.4},
$\dim\{S^2_{\mathcal{U}_\pm^\star}(\mathfrak{H}_\pm)\}\le 4$. Theorem \ref{thm-1.4} now follows from
Lemma~\ref{lem-2.1} and from Lemma~\ref{lem-2.3}.
\hfill\qed

\medbreak\noindent{\bf Proof of Theorem \ref{thm-1.1}.}
Let $(M,g,J_\pm)$ be an almost para/pseudo-Hermitian manifold of
dimension $m\ge6$ (the case $m=4$ is analogous). We consider variations
$(M,g,J_\pm^{\varTheta})$. Subtracting
$\nabla\Omega_\pm(M,g,J_\pm)(P)$ has no effect on the question of surjectivity. Every $\vartheta\in\mathfrak{o}\otimes T^*M$
can be written in the form $\vartheta=d\varTheta(P)$ for some admissible $\varTheta$. Thus it suffices to show
$\Xi_\pm(\mathfrak{o})=\mathfrak{H}_\pm$. By Lemma \ref{lem-5.1},
$\Xi(\mathfrak{o})$ is not perpendicular to
$W_{\pm,i}$ for $1\le i\le 4$. By Theorem \ref{thm-1.4}, $W_{\pm,i}$ is an irreducible submodule of $\mathfrak{H}_\pm$ which
occurs with multiplicity $1$. Thus by Lemma \ref{lem-2.1},
$W_{\pm,i}\subset\Xi(\mathfrak{o})$ for
$1\le i\le 4$. Theorem \ref{thm-1.4} now shows $\mathfrak{H}_\pm\subset\Xi(\mathfrak{o})$ as desired.
\hfill\qed

\section{Varying the metric}\label{sect-6}
Let $(M,g,J_\pm)$ be a para/pseudo-Hermitian manifold. Fix $P$ in $M$ and let
$$(V,\langle\cdot,\cdot\rangle,J_\pm^0):=(T_PM,g_P,J_{\pm,P})\,.$$
Let $\mathfrak{gl}_\pm$ be the Lie algebra of
$\operatorname{GL}_\pm$ at $P$. Given $\tilde\vartheta\in\mathfrak{gl}\otimes V^*$, we may find a smooth
map
$\tilde\varTheta$ from a neighborhood of $P$ in $M$ to $\operatorname{GL}_\pm$ so that $\tilde\varTheta(P)=\id$, so that
$\tilde\varTheta=\id$ away from a neighborhood of $P$, and so that $d\tilde\varTheta(P)=\tilde\vartheta$.
We define a new pseudo-Riemannian metric $g^{\tilde\varTheta}$ which agrees with $g$ at $P$ and which agrees with $g$
away from a neighborhood of $P$ by setting:
$$g^{\tilde\varTheta}(x,y)=(\tilde\varTheta x,\tilde\varTheta y)\,.$$
Since $\varTheta J_\pm=J_\pm\varTheta$, $g^{\tilde\varTheta}$ is a para/pseudo-Hermitian metric. Set:
$$
\tilde\Xi_\pm(\tilde\vartheta):=\left\{\nabla\Omega_\pm(V,g^{\tilde\varTheta},J_\pm)-\nabla\Omega_\pm(V,g,J_\pm^0)\right\}(P)\,.
$$
We may then use Lemma \ref{lem-3.2} to see that $\tilde\Xi_\pm(\tilde\vartheta)\in W_{\pm,3}$ is independent of the choice of
$\tilde\varTheta$ and defines a $\mathcal{U}_\pm^\star$ module morphism from $\mathfrak{gl}_\pm\otimes V^*\otimes\chi$ to
$\mathfrak{H}_\pm$ by computing:
\begin{eqnarray*}
&&\tilde\Xi_\pm(\tilde\vartheta)(x,y;z)\\
&=&\textstyle\frac12\left\{\langle\tilde\vartheta(J_\pm y)x,z\rangle+\langle x,\tilde\vartheta(J_\pm y)z\rangle
+\langle\tilde\vartheta(y)J_\pm x,z\rangle+\langle J_\pm x,\tilde\vartheta(y)z\rangle\right.\\
&&\quad\left.-\langle\tilde\vartheta(J_\pm x)y,z\rangle+\langle y,\tilde\vartheta(J_\pm x)z\rangle
+\langle\tilde\vartheta(x)J_\pm y,z\rangle+\langle J_\pm y,\tilde\vartheta(x)z\rangle\right\}\,.
\end{eqnarray*}

Thus to prove Theorem \ref{thm-1.2}, it suffices to show that $\Xi_\pm$ is surjective. Since we have subtracted the effect of
the background metric, we may take the flat metric $g=\langle\cdot,\cdot\rangle$. As in Section
\ref{sect-5}, we introduce a normalized orthonormal basis $\{e_1,\dots,e_{\bar m},f_1,\dots,f_{\bar m}\}$ for $V$. Let $\alpha$
be a smooth function which is compactly supported near $P=0$ with $\alpha(0)=0$ and $d\alpha(0)=dx^1$. Set:
$$\tilde\varTheta e_i=\left\{\begin{array}{rll}
e^\alpha e_i&\text{if}&i=1,2\\ e_i&\text{if}&i\ge3\end{array}\right\}\quad\text{and}\quad
\tilde\varTheta f_i=\left\{\begin{array}{rll}
e^\alpha f_i&\text{if}&i=1,2\\ f_i&\text{if}&i\ge3\end{array}\right\}\,.
$$
Let $\tilde\vartheta=d\tilde\varTheta(0)=\tilde\vartheta_0\otimes dx^1$ where $\tilde\vartheta_0$ is orthogonal projection on
$\operatorname{Span}\{e_1,e_2,f_1,f_2\}$:
\begin{eqnarray*}
\tilde\vartheta_0e_i&=&\left\{\begin{array}{lll}e_i&\text{if}&i=1,2\\0&\text{if}&i\ge3\end{array}\right\}\quad\text{and}\quad
\tilde\vartheta_0f_i=\left\{\begin{array}{lll}f_i&\text{if}&i=1,2\\0&\text{if}&i\ge3\end{array}\right\} .
\end{eqnarray*}
The associated metric takes the form:
\begin{eqnarray*}
g^{\tilde\varTheta}_\pm&=&
e^{2\alpha}\varepsilon_1(e^1\otimes e^1\mp f^1\otimes f^1)+e^{2\alpha}\varepsilon_2(e^2\otimes e^2\mp f^2\otimes f^2)\\
&+&\sum_{i\ge3}\varepsilon_i(e^i\otimes e^i\mp f^i\otimes f^i)\,.
\end{eqnarray*}
Set $H_\pm := \nabla\Omega_\pm(0) = \tilde\Xi_\pm(\vartheta)$. We use Lemma \ref{lem-3.1} to see $\tau_1(H_\pm) = 2e^1$ and
thus $H_\pm$ has a non-trivial component in $W_{\pm,4}$. Since $H_\pm(e_1, e_3; f_3) = 0$ and
$\sigma_\pm(e^1)(e_1, e_3; f_3) \neq 0$, $H_\pm$ also has a non-zero component in $W_{3,\pm}$. Theorem \ref{thm-1.2} now follows.\hfill\qed

\section{The 16 classes of almost pseudo-Hermitian manifolds}\label{sect-7}

\medbreak\noindent{\bf Proof of Theorem \ref{thm-1.6}.}
If $(M,g,J_-)$ is a $\xi$-manifold, then $(M,-g,J_-)$ also is a $\xi$-manifold.
Thus by replacing
$g$ by
$-g$ if need be, we may assume without loss of generality that $p\le q$ and consequently, as $m\ge10$, that $6\le q$ to
establish Theorem \ref{thm-1.6}.
We shall use product structures. The projections $\pi_{i,-}$ for $i=1,2,3$ and the map $\tau_1$
are compatible with Cartesian product; the splitting $\sigma_-$ is not. This causes a small amount of additional
technical fuss.

Suppose first that $W_4\not\subset\xi$. By Theorem \ref{thm-1.5} we may choose a $\xi$-manifold  $(M_1,g_1,J_{1,-})$ of
Riemannian signature $(0,q)$. Let
$(M_2,g_2,J_{2,-})$ be a flat K\"ahler torus of signature $(p,0)$. Let
\begin{equation}\label{eqn-7.a}
M=M_1\times\mathbb{T}^{(p,0)},\qquad g:=g_1+g_2,\quad J_-=J_{1,-}\oplus J_{2,-}\,.
\end{equation}
Then $(M,g,J_-)$ is an almost pseudo-Hermitian manifold of signature $(p,q)$. We have
$\nabla\Omega_g=\nabla\Omega_{g_1}$ and $\tau_1(\nabla\Omega_g)=\tau_1(\nabla\Omega_{g_1})=0$. Thus
$\pi_{3,-}\nabla\Omega_g$ is projection on $W_{-,3}$; this would not be the case if $\tau_1$
was non-zero and this fact played an important role in the analysis of Section \ref{sect-6}. Since
$\pi_{i,-}\nabla\Omega_g=\pi_{i,-}\nabla\Omega_{g_1}$, it now follows that
$(M,g,J_-)$ is a
$\xi$ manifold in this special case.

Next we suppose that $\xi=\eta\oplus W_{-,4}$. Let $(M,g,J_-)$ be an $\eta$-manifold of signature $(p,q)$. We make a
conformal change of metric and set
$\tilde g:=e^{2f}g$; it then follows from Lemma \ref{lem-3.2} that
$$\nabla\Omega_{\tilde g}=e^{2f}\nabla\Omega_g-e^{2f}\sigma_{-,g}(df)$$
where we use the original metric to define the splitting $\sigma_{-,g}$. This has
a non-trivial $W_{4,-}$ component and the components $W_{i,-}$ for $1\le i\le 3$ are not affected.
\hfill\qed

\section*{Acknowledgments}
Research of the authors partially supported by project MTM2009-07756 (Spain) and by INCITE09 207 151 PR (Spain).


\begin{thebibliography}{aaa}

\bibitem{AGK11}
Alekseevsky D.,  Guilfoyle B., and Klingenberg W.,
``On the Geometry of Spaces of Oriented Geodesics",
arXiv:0911.2602v1 [math.DG].


\bibitem{AMT09} Alekseevsky D., Medori C., and Tomassini A.,
``Homogeneous para-Kaehler Einstein manifolds",
{\it Russian Math. Surveys} {\bf 64} (2009), 1--43.



\bibitem{AFS05} Alexandrov B., Friedrich Th., and Schoemann N.,
``Almost Hermitian $6$-manifolds revisited",
{\it J. Geom. Phys.} \textbf{53} (2005), 1--30.

\bibitem{Arm02}
Armstrong, J.,
``An ansatz for almost-K\"{a}hler, Einstein $4$-manifolds",
{\it J. Reine Angew. Math.}  {\bf 542}  (2002), 53?84.

\bibitem{BLS07}
Bor G., Hern\'andez-Lamoneda L., and Salvai M.,
``Orthogonal almost-complex structures of minimal energy",
{\it Geom. Dedicata} {\bf 127} (2007), 75--85.


\bibitem{BGN11} Brozos-V\'azquez M., Gilkey P., and Nik\v cevi\'c S., ``Geometric Realizations of Curvature",
Imperial College Press (2011).



\bibitem{CI08}
Cleyton, R. and Ivanov, S.,
``Conformal equivalence between certain geometries in dimension 6 and 7",
{\it Math. Res. Lett.}  {\bf 15}  (2008),   631-?640.

\bibitem{CMMS04} Cort\'{e}s V.,  Mayer Ch.,  Mohaupt Th., and Saueressig F.,
``Special geometry of Euclidean supersymmetry 1. Vector multiplets",
{\it J. High Energy Phys.} \textbf{03} (2004) 028.

\bibitem{CMMS05}  Cort\'{e}s V., Mayer Ch.,  Mohaupt Th., and  Saueressig F.,
`` Special geometry of Euclidean supersymmetry II. Hypermultiplets and the c-map",
{\it J. High Energy Phys.} \textbf{06} (2005) 025.

\bibitem{CS07}
Cort\'es V. and Sch\"afer L.
``Flat nearly K\"ahler manifolds'',
{\it Ann. Global Anal. Geom.} {\bf 32} (2007), 379--389.


\bibitem{EPS09}
Euh Y., Park J. H., and Sekigawa K.,
``Nearly K\"ahler manifolds with vanishing Tricerri-Vanhecke Bochner curvature tensor",
{\it Differential Geom. Appl.} {\bf 27} (2009), 250--256.

\bibitem{FFS94} Falcitelli M., Farinola A., and Salamon S.,
``Almost-Hermitian geometry",
{\it Differential Geom. Appl.} \textbf{4} (1994), 259--282.


\bibitem{Fu58} Fukami T., ``Invariant tensors under the real representation of unitary groups
and their applications", {\it J. Math. Soc. Japan} {\bf 10} (1958), 135--144.


\bibitem{GM91} Gadea P. and Masque J.,
``Classification of almost para-Hermitian manifolds",
{\it Rend. Mat. Appl.} {\bf 11} (1991), 377--396.

\bibitem{GM00}  Garc\'{\i}a-R\'{\i}o E. and  Matsushita Y.,
``Isotropic K\"{a}hler structures on Engel $4$-manifolds'',
\emph{J. Geom. Phys.}  \textbf{33}  (2000),   288--294.

\bibitem{GPSR07}
Gran U., Papadopoulos G., Sloane P., and Roest D.,
``Geometry of all supersymmetric type I backgrounds'',
{\it J. High Energy Phys.} (2007) 074.

\bibitem{GH80} Gray A. and Hervella L.,
``The sixteen classes of almost Hermitian manifolds and their linear invariants",
{\it Ann. Mat. Pura Appl.} {\bf 123} (1980), 35--58.

\bibitem{GK05}  Guilfoyle B. and  Klingenberg W.,
``An indefinite K\"{a}hler metric on the space of oriented lines",
\emph{J. London Math. Soc.} \textbf{72} (2005), 497--509.


\bibitem{HL00}
Hern\'{a}ndez-Lamoneda, L.,
``Curvature vs. almost Hermitian structures",
{\it Geom. Dedicata} {\bf 79}  (2000),  205?218.

\bibitem{Iw58}
Iwahori N., ``Some remarks on tensor invariants of $O(n)$, $U(n)$, $Sp(n)$",
{\it J. Math. Soc. Japan} {\bf 10} (1958), 146--160.



\bibitem{MS04} Mart\'{\i}n-Cabrera F. and  Swann A.,
``Almost Hermitian structures and quaternionic geometries",
{\it Differential Geom. Appl.} \textbf{21} (2004), 199--214.

\bibitem{MHL07}  Matsushita Y.,  Haze S., and  Law P. R.,
``Almost K\"{a}hler-Einstein structures on $8$-dimensional Walker manifolds",
\emph{Monatsh. Math.} \textbf{150} (2007), 41--48.

\bibitem{MO08}
Moroianu A. and Ornea L.,
``Conformally Einstein products and nearly K\"ahler manifolds",
{\it Ann. Global Anal. Geom.} {\bf 33} (2008), 11--18.

\bibitem{Na83}
Naveira, A. M.;
``A classification of Riemannian almost-product manifolds",
{\it Rend. Mat. (7)} {\bf 3}  (1983),  577-?592.

\bibitem{P08}
Papadopoulos G.,
``Killing-Yano equations and G structures'',
{\it Classical Quantum Gravity} {\bf 25} (2008), no. 10, 105016, 8 pp.

\bibitem{MN03}  San Martin L. A. B. and Negreiros C. J. C.,
``Invariant almost Hermitian structures on flag manifolds",
{\it Adv. Math.} \textbf{178} (2003), 277--310.


\bibitem{SY08}  Sekigawa K. and  Yamada A.,
``Compact indefinite almost K\"{a}hler Einstein manifolds",
\emph{Geom. Dedicata} \textbf{132} (2008), 65--79.



\bibitem{TV81} Tricerri F. and Vanhecke L.,
``Curvature tensors on almost Hermitian manifolds", {\it  Trans. Amer. Math. Soc.}
{\bf  267}  (1981), 365--397.

\bibitem{TV83} Tricerri F. and Vanhecke L.,
``Homogeneous structures on Riemannian manifolds", London Math. Soc. Lecture Note Ser., 83, Cambridge University Press, Cambridge, 1983.



\bibitem{W46} Weyl H., ``The classical groups", Princeton Univ. Press,
Princeton (1946) (8${}^{\operatorname{th}}$ printing).



\end{thebibliography}
\end{document}